\documentclass[a4paper,12pt]{article}
\usepackage{amsfonts,amsmath,amsxtra}

\newtheorem{theorem}{Theorem}
\newtheorem{corollary}{Corollary}

\newenvironment{definition}
{\smallskip\noindent{\bf Definition\/}:}{\smallskip\par}
\newenvironment{proposition}
{\smallskip\noindent{\bf Proposition\/}.}{\smallskip\par}
\newenvironment{example}
{\smallskip\noindent{\bf Example\/}.}{\smallskip\par}

\newenvironment{proof}
{\noindent{\bf Proof\/}.}{{ $\Box$}\smallskip\par}

\title{Adams operations and power structures.}
\author{E. Gorsky \footnote{Partially supported by the grants RFBR-007-00593, INTAS-05-7805, NSh-4719.2006.1, and the Moebius Contest fellowship for the young scientists.}}
\date{}

\begin{document}

\maketitle

\section{Introduction}

We construct a family of additive endomorphisms $\Psi_k, k=1, 2\ldots$ of the Grothendieck ring of quasiprojective varieties and the Grothendieck ring of Chow motives. Each of these maps being applied to a given variety gives a polynomial of its symmetric powers with integer coefficients. For example, $\Psi_1(X)=X, \Psi_2(X)=2[S^{2}X]-[X^2],$ where $S^2X$ is a symmetric square of $X$. For a polynomial of the affine line there is a formula
$$\Psi_k(P(\mathbb{L}))=P(\mathbb{L}^k).$$
The construction of these maps has a lot in common with the construction of the Adams operations in the K-theory,
and for the Grothendieck ring of Chow motives two additional sets of equations
$$\Psi_i(X\cdot Y)=\Psi_i(X)\cdot\Psi_i(Y), \Psi_i\circ\Psi_j=\Psi_{ij},$$
analogous to the ones for the Adams operations, hold. This fact follows from the specialty of the 
 $\lambda$-structure over the Grothendieck ring of motives proved by F. Heinloth (\cite{hein}).

These operations are used for study of the so-called power structure over the Grothendieck ring, constructed by S. Gusein-Zade, I. Luengo and A. Melle-Hernandez (\cite{powers}). We prove the inversion formula, which provides a possibility to express  explicitly the exponents $B_i$ via the coefficients $A_j$ in a formula
$$1+A_1\cdot t+A_2\cdot t^2+\ldots=(1-t)^{-B_1}(1-t^2)^{-B_2}\cdots,$$
where the right hand side is considered in a sense of a power structure.
This formula looks very similar to the "motivic Moebius inversion" of D. Bourqui (\cite{bourqui}).

As an example we calculate the class of the variety of irreducible polynomials of a given degree of arbitrary number of variables. This gives a way to compute, for example, all Hodge-Deligne numbers of this variety.

Moreover, we give a clear interpretation on a language of power structures of the E. Getzler's formulas
(\cite{getzler}) for the characters of the equivariant cohomologies of the natural symmetric group action on the configuration spaces of ordered tuples of points on a given variety.

We also recall some known results concerning power structures like the alternative proof of the L. Goettsche's formula for the Betti numbers of the Hilbert schemes of points on a surface from \cite{hilbert}.

\section{Power structures}

The notion of a power structure over a (semi)ring was introduced by S. Gusein-Zade, I. Luengo and A. Melle-Hernandez in \cite{powers}.

\begin{definition}
A power structure over a ring $R$
is a map $$(1+tR[[t]])\times R\rightarrow 1+tR[[t]]:(A(t),m)\mapsto (A(t))^{m},$$
satisfying the following equations:

1.$(A(t))^{0}=1,$

2.$(A(t))^{1}=A(t),$

3.$((A(t)\cdot B(t))^{m}=((A(t))^{m}\cdot((B(t))^{m},$

4.$(A(t))^{m+n}=(A(t))^{m}\cdot(A(t))^{n},$

5.$(A(t))^{mn}=((A(t))^{n})^{m},$

6.$(1+t)^{m}=1+mt+$ terms of the higher degree,

7.$(A(t^{k}))^{m}=((A(t))^{m})|_{t\rightarrow t^{k}}.$
\end{definition}

A power structure is said to be finitely determined,
if for every 

$N>0$ there exists $M>0$, such that
$N$ -- jet of a series $(A(t))^{m}$ is uniquely determined by the
$M$ -- jet of the series $A(t)$.

In \cite{powers} it is proved that a finitely determined power structure is defined if and only if one has a rule
defining $(1-t)^{-m}$
for every $m\in R$, and $$(1-t)^{-m-n}=(1-t)^{-m}\cdot(1-t)^{-n}.$$

If we have a series $1+A_1t+A_2t^2+A_3t^3+\ldots$,
then, after dividing it by $(1-t)^{-A_1}$, we'll get a series of a form $1+D_2t^2+D_3t^3+\ldots$, then we divide it by $(1-t^2)^{-D_2}$ etc.

Finally we'll get a decomposition of a power series into an infinite product
\begin{equation}
\label{prod}
1+A_1\cdot t+A_2\cdot t^2+\ldots=(1-t)^{-B_1}(1-t^2)^{-B_2}\cdots.
\end{equation}
If all series of a form $(1-t^{k})^{-B_k\cdot X}$ are known, it's easy to compute $A(t)^X$ by setting
$$(1+A_1\cdot t+A_2\cdot t^2+\ldots)^X=\prod_{k=1}^{\infty}(1-t^k)^{-B_k\cdot X}.$$

\medskip

If a ring $R$ is a $\mathbb{Q}$-algebra, one can define the exponential and the logarithmic maps, so one can define a "usual" power structure by the formula
$$(1+A_1\cdot t+A_2\cdot t^2+\ldots)^X:=\exp(X\cdot\ln(1+A_1\cdot t+A_2\cdot t^2+\ldots)).$$
It's easy to see that the equations 1-7 are satisfied.

It is important to note that this power structure is not unique. Below we'll discuss a bunch of examples of useful and important power structures which are far from this one. Most of their definitions do not use a division by integers, so
they are defined over $\mathbb{Z}$, not only for $\mathbb{Q}$-algebras.

\subsection{Power structure over the Grothendieck ring of varieties}

By $K_0(Var_{\mathbb{C}})$ we denote the Grothendieck ring of 
quasiprojective algebraic varieties. It is generated by the isomorphism
classes of complex quasiprojective algebraic varieties modulo relations of form
 $[X]=[Y]+[X\setminus Y]$ where $Y$ is a Zariski closed subset of $X$.
Multiplication is given by the formula $[X]\cdot [Y]=[X\times Y].$
Let $\mathbb{L}\in K_0(Var_{\mathbb{C}})$ denote the class of the affine line.

Over the Grothendieck ring of quasiprojective varieties
one can define a power structure by the formula
$$(1-t)^{-[X]}=1+[S^{1}X]t+[S^{2}X]t^{2}+\ldots,$$
where $S^{k}X$ denotes the $k$th symmetric power of $X$.
Analogously one can define a power structure over the Grothendieck ring of Chow motives.
For example, for
$j\ge 0$ it is known (\cite{powers}) that $[S^{k}\mathbb{L}^{j}]=\mathbb{L}^{kj}$, so
$$(1-t)^{-\mathbb{L}^j}=\sum_{k=0}^{\infty}t^k\mathbb{L}^{kj}=(1-t\mathbb{L}^j)^{-1}.$$

\medskip

This power structure has a nice and clear geometric meaning: if $A_1, A_2,\ldots$ and $X$  are some quasiprojective varieties and $$(1+A_1\cdot t+A_2\cdot t^2+\ldots)^X=1+B_1\cdot t+B_2\cdot t^2+\ldots,$$
then the varieties $B_n$ have a following geometric description. Consider a function $I$ on a disjoint union $\sqcup_{i=1}^{\infty}A_{i}$ which equals identically to $i$ on $A_{i}$. Then $B_n$ is a set of pairs $(K,\varphi)$, where $K$ is a finite subset of $X$ and $\varphi:K\rightarrow \sqcup_{i=1}^{\infty}A_i$ is a map such that $$\sum_{x\in K}I(\varphi(x))=n.$$
This set of pairs can be naturally  equipped with a structure of a quasiprojective algebraic variety.
Furthermore, from this geometric definition it is easy to check that this construction satisfies all properties of a power structure.

Less formally this construction can be described in a following way: on a variety $X$ there live particles equipped with  some natural numbers (multiplicities, masses, charges ...). A particle of a given multiplicity $k$ has a complicated space of internal states which are parametrized by points of a quasiprojective variety $A_k$. Now $B_n$ is 
a configuration space of tuples of particles of total multiplicity $n$. For example, if all $A_k$ are one-point sets, then $B_n$ consists of all possible tuples of distinct points on $X$ with multiplicities of total multiplicity $n$, that is $B_n=S^{n}X$.
Hence
$$1+[S^{1}X]t+[S^{2}X]t^{2}+\ldots=(1+t+t^2+\ldots)^{[X]}=((1-t)^{-1})^{[X]}=(1-t)^{-[X]}.$$
If $A_1$ is a point, and all $A_i$ are empty, then $B_n$ is a set of unordered tuples of distinct points on $X$. 
A generating function for the classes of these sets has a form $(1+t)^{[X]}$.

Less trivial examples of usage of power structures are also known. For example, the Jordan normal form of a matrix is a 
tuple of its eigenvalues with attached partitions (Young diagrams).
Therefore the set of Jordan forms of $n\times n$ matrices is a coefficient with number $n$ in the series $$(1+A_1\cdot t+A_2 t^2+\ldots)^{X},$$
where $X$ is a set consisting of all possible eigenvalues (whole $\mathbb{C}$ or $\mathbb{C}^{*}$, if we consider only nondegenerate matrices), and $A_k$ is a set of Young diagrams of weight $k$. This example produce some curious combinatorial identities which can be found  (in a slightly different terminology) in the article (\cite{rodrigues}).

Another application of the technique of the power structures is the geometry and the combinatorics of the Hilbert schemes of points on varieties. If $X$ is a smooth projective variety of dimension $d$, then the set of its zero-dimensional subschemes of length $k$ can be equipped with a structure of a projective variety. It is said to be a Hilbert scheme of $k$ points on $X$ and is denoted by $Hilb^{k}(X)$. It turns out that the following identity is true (\cite{hilbert}):
$$(1+[Hilb^{1}(\mathbb{C}^d,0)]\cdot t+[Hilb^{2}(\mathbb{C}^d,0)]\cdot t^2+\ldots)^{[X]}=1+[Hilb^{1}(X)]\cdot t+[Hilb^{2}(X)]\cdot t^2+\ldots,$$
where $Hilb^{k}(\mathbb{C}^d,0)$ is a Hilbert scheme parametrizing subschemes of  $\mathbb{C}^d$ of length $k$ 
with the support at the origin. A zero-dimensional subscheme of $X$ can be considered as a pair $(K,\varphi)$, where $K$ is a finite subset of $X$, and $\varphi$ is a map from $K$ to the set of zero-dimensional subsets of $X$ with one-point support, what motivates this formula. 
This identity does not follow directly from the definition of the power structure, but it can be deduced from it after some technical work (\cite{hilbert}).

Identities like this holds in the Grothendieck ring of varieties, so they look quite abstract. Nevetherless, they have concrete and very powerful geometric corollaries. For example, the Euler characteristic is an additive invariant of algebraic varieties in a sense that
 $$\chi(X)=\chi(Y)+\chi(X\setminus Y),$$
if $Y$ is a Zariski closed subset of $X$. Moreover, $\chi(X\times Y)=\chi(X)\cdot\chi(Y),$
so $\chi:K_0(Var_{\mathbb{C}})\rightarrow \mathbb{Z}$ is a ring homomorphism.
There exist some other invariants of algebraic varieties. For example, the Hodge-Deligne polynomial (\cite{dankho})
is a ring homomorphism
$$e:K_0(Var_{\mathbb{C}})\rightarrow \mathbb{Z}(u,v),$$
which coincides for smooth projective varieties with the generating function for the Hodge numbers:

$$e(X)=\sum_{i,j}(-1)^{i+j}u^{i}v^{j}h_{i,j}(X).$$

Over the polynomial ring there exists a power structure defined by the equation
$$(1-t)^{-\sum a_{\underline{k}}\underline{x}^{\underline{k}}}=
\prod (1-t\underline{x}^{\underline{k}})^{-a_{\underline{k}}}.$$

In \cite{powers} it is proved that the Hodge-Deligne polynomial is a morphism of power structures, that is, if $M,A_1,A_2\ldots\in K_0(Var_{\mathbb{C}})$, then
$$e((1+A_1t+A_2t^2+\ldots)^M)=(e(1+A_1t+A_2t^2+\ldots))^{e(M)}.$$
In particular, $$\chi((1+A_1t+A_2t^2+\ldots)^M)=(1+\chi(A_1)t+\chi(A_2)t^2+\ldots))^{\chi(M)}.$$

Therefore we can translate equations in the Grothendieck ring into equations involving Hodge numbers (or Euler characteristics) of some concrete varieties. For example, in the example with Hilbert schemes of points from the identity in the Grothendieck ring one can obtain the formula of L. Goettsche
(\cite{goettsche}) for the Betti numbers of the Hilbert scheme of points on a surface. In the notation of the power structure over the polynomial ring it has a form :
$$(1+e(Hilb^{1}(\mathbb{C}^d,0))\cdot t+e(Hilb^{2}(\mathbb{C}^d,0))\cdot t^2+\ldots)^{e(X)}=1+e(Hilb^{1}(X))\cdot t+e(Hilb^{2}(X))\cdot t^2+\ldots.$$

\subsection{Special $\lambda$ -- rings}

Let $R$ be a commutative ring with a unity.

\begin{definition}
A map $\lambda_t:R\rightarrow 1+tR[[t]]$ is said to be a $\lambda$-structure over $R$, if $$\lambda_t(X+Y)=\lambda_t(X)\lambda_t(Y)$$
and $\lambda_t(X)=1+tX+\ldots.$
\end{definition}

If a  $\lambda$-structure over a ring is given, one can define a power structure over it, setting  $(1-t)^{-X}=\lambda_t(X).$
From the other hand, a power structure over a ring induces a bunch of different power structures of a form
$$\lambda^{A}_{t}(X)=(A(t))^{X},$$
where $A(t)$ is an arbitrary series of a form $1+t+\ldots$.
Nevetherless, below, if the contrary is not said, a $\lambda$-structure over a ring with a power structure 
will be considered equal to $(1-t)^{-X}$.

Let $\sigma_{i}$ be an $i$th elementary symmetric polynomial of variables  $\xi_1,\ldots,\xi_N$, $s_i$  is $i$th elementary symmetric polynomial of variables $x_1,\ldots,x_N$.
Let 

$P_n(\sigma_1,\ldots,\sigma_n;s_1,\ldots,s_n)$ be a coefficient at $t^n$ in the series $$\prod_{1\le i,j\le N}(1+\xi_i x_j t),$$ and $P_{n,r}(\sigma_1,\ldots, \sigma_{nr})$ is a coefficient at
$t^n$ in the series $$\prod_{1\le i_1<\ldots<i_r\le N}(1+\xi_{i_1}\ldots\xi_{i_r}t).$$

Consider a following $\lambda$-structure over a ring $1+tR[[t]]$.
Addition  is given by multiplication, multiplication  $\circ$ is given by a formula
$$(1+a_1t+a_2t^2+\ldots)\circ(1+b_1t+b_2t^2+\ldots)=1+\sum_{n=1}^{\infty}P_n(a_1,\ldots, a_n;b_1,\ldots,b_n)t^n,$$
and $\lambda$-structure is given by the formula
$$\Lambda_r(1+a_1t+a_2t^2+\ldots)=1+\sum_{n=1}^{\infty}P_{n,r}(a_1,\ldots,a_{nr})t^n.$$

\begin{definition}
$\lambda$-structure over a ring $R$ is said to be {\bf special}, if
$\lambda_t:R\rightarrow 1+tR[[t]]$ is a ring homomorphism, preserving the $\lambda$-structure.
\end{definition}

F. Heinloth (\cite{hein}) proved that the $\lambda$-structure over the Grothendieck ring of Chow motives is special. 
It is not known if the analogous statement holds for the Grothendieck ring of varieties. Since the natural additive maps like the Hodge-Deligne polynomial factorizes through the Chow motives, it's not a big problem.

To use the theorem of Heinloth for calculations, it is useful to reformulate, following
 \cite{knut}, the definition of a special $\lambda$-structure.

Let $$\lambda_t(X)^{-1}\cdot {d\over dt}\lambda_t(X)=\sum_{n=1}^{\infty}\Psi_n(X)t^{n-1}.$$

The definition of a  $\lambda$-structure is equivalent to the identity $$\Psi_i(X+Y)=\Psi_i(X)+\Psi_i(Y),$$ and if a $\lambda$-structure is special, then (\cite{knut}) $$\Psi_i(XY)=\Psi_i(X)\Psi_i(Y)\,\,\,\, \mbox{\rm and}\,\,\,\, \Psi_i\circ\Psi_j=\Psi_{ij}$$
for all $i$ and $j$.

\begin{definition}
These homomorphisms $\Psi_{k}$ are said to be Adams operations on a ring $R$.
\end{definition}

\begin{example}
Let $X=\mathbb{L}^{k}.$ Then $\lambda_t(X)=(1-\mathbb{L}^{k}t)^{-1},$
so $$\lambda_t(X)^{-1}\cdot {d\over dt}\lambda_t(X)={\mathbb{L}^{k}\over 1-\mathbb{L}^{k}t},$$
hence $$\Psi_{n}(\mathbb{L}^{k})=\mathbb{L}^{kn}.$$
Since $\Psi_n$ are additive operations,  $\Psi_n(P(\mathbb{L}))=P(\mathbb{L}^{n})$, therefore, 
over the subring of polynomials of $\mathbb{L}$ the $\lambda$-structure is special.
\end{example}

In fact, the Adams operation on the polynomial ring of arbitrary number of variables have a similar form.
If $q(\underline{x})=\sum a_{\underline{k}}\underline{x}^{\underline{k}}$, then
$$(1-t)^{-q(\underline{x})}=(1-t)^{-\sum a_{\underline{k}}\underline{x}^{\underline{k}}}=
\prod (1-t\underline{x}^{\underline{k}})^{-a_{\underline{k}}},$$
hence
$$(1-t)^{q(\underline{x})}{d\over dt}(1-t)^{-q(\underline{x})}=\sum {a_{\underline{k}}x^{\underline{k}}\over 1-t\underline{x}^{\underline{k}}}=
\sum_{n=1}^{\infty}q(x_1^n, x_2^n,\ldots)t^{n-1},$$
so $$\Psi_{n}(q(x_1,x_2,\ldots))=q(x_1^n,x_2^n,\ldots).$$
For example, this means that this $\lambda$-structure over the polynomial ring is special.

\begin{example}
Another important example of a ring with a natural special
$\lambda$-structure is a Grothendieck ring of representations of a given finite group $G$. Each representation
$\rho:G\rightarrow GL(V)$ is in the 1-to-1 correspondence with its character
$\chi(g)=\mbox{\rm tr}\,\,\,\rho(g),$ and the Grothendieck ring of representations is isomorphic to a ring of functions on $G$ invariant under the conjugation. The character of the sum (tensor product) of two representations is equal to a sum (product) of the characters of these representations, so the character map is a ring homomorphism.

Over the ring of representations we have a natural $\lambda$-structure:
$$\lambda_t(V)=1+\sum_{k=1}^{\infty}S^{k}V\cdot t^k.$$
Let $\xi_1,\ldots,\xi_n$ are eigenvalues of $\rho(g)$ for some $g\in G$.
The eigenvalues of the operator $S^{k}\rho(g)$ acting on the space $S^{k}V$ are all products of a form
$\xi_{i_1}\ldots\xi_{i_k}, i_1\le\ldots\le i_k$. Therefore,
$$1+\sum_{k=1}^{\infty}t^k\mbox{\rm tr}S^{k}\rho(g)=\prod_{j=1}^{n}(1-\xi_jt)^{-1},$$
so the value of the character of $\Psi_k(\rho)$ at the element $g$ is equal to
$\xi_1^{k}+\ldots+\xi_n^{k}=\chi_{\rho}(g^{k}).$
Therefore $$\Psi_{k}\chi(g)=\chi(g^k),$$
and, for example, the power structure is special.
\end{example}

\begin{example}
First appearance of the Adams operations was in the K-theory (\cite{adams}). Let $X$ be an arbitrary topological space, $K_0(X)$ is a Grothendieck group of vector bundles over it. Over $K_0(X)$ there is a natural  $\lambda$-structure: if $E$ is a (virtual) bundle, then
$$(1-t)^{-E}=1+\sum_{k=1}^{\infty}S^{k}E\cdot t^k.$$
For the calculation of the Adams operations one can use the decomposition principle: if $E=\oplus_{i=1}^{n}E_i,\,\, \mbox{\rm rk}\,\, E_i=1$, then
$$\Psi_k(E)=\oplus_{i=1}^{n}E_i^{k}.$$
The proof of this fact is completely analogous to the previous example, and the power structure is special.
\end{example}

\begin{example}
Let $R=\oplus_{i=0}^{\infty}R_i$ be a graded ring,
$R_i\cdot R_j\subset R_{i+j}$. For every  $x\in R_j$ let
$$\Psi_k(x)=k^j\cdot x.$$
It is easy to see that these operations are ring homomorphisms and $$\Psi_k(\Psi_m(x))=\Psi_{km}(x).$$ Therefore $\Psi_k$
are Adams operations for some special power structure. It is clear how to reconstruct this structure: if $x\in R_j$,
then $$(1-t)^{-x}=\exp(\sum_{k=1}^{\infty}\Psi_k(x){t^k\over k})=
\exp(\sum_{k=1}^{\infty}k^{j-1}t^k)=1+t+(2^{j-1}+{1\over 2})t^2+(3^{j-1}+2^{j-1}+{1\over 6})t^3+\ldots.$$

This structure is strange at a first glance, but it appears naturally, for example, in the even-dimensional cohomologies
of an arbitrary topological space $X$. The Chern character $ch$ is a homomorphism from $K_0(X)$ to $H^{2*}(X)$. If $E$ is a line bundle over $X$, then $$c_1(\Psi_k(E))=c_1(E^k)=kc_1(E)=\Psi_k(c_1(E)).$$ Since $\Psi_k$ are ring homomorphisms,
$$ch(\Psi_k(E))=e^{c_1(\Psi_k(E))}=\Psi_k(e^{c_1(E)})=\Psi_k(ch(E)).$$ From the decomposition principle and the properties of the Adams operations it follows that for every bundle $E$ $$ch(\Psi_k(E))=\Psi_k(ch(E)),\,\,\,ch((1-t)^{-E})=(1-t)^{-ch(E)}.$$
\end{example}

\section{An inversion formula}

The proof of the formula (\ref{prod}), cited from \cite{powers}, is clear, but it does not give any explicit formulas  expressing $B_i$ through the coefficients of  $A(t)$. It turns out that such formulas can be written in terms of the Adams operations.

\begin{theorem}
Let $A(t)^{-1}{d\over dt}A(t)=\sum_{n=1}^{\infty}C_{n}t^{n-1}.$
Then $$nB_n=\sum_{i_1,\ldots,i_{s-1}>1,i_1\cdot\ldots\cdot i_s=n}(-1)^{s-1}\Psi_{i_1}\ldots\Psi_{i_{s-1}}(C_{i_s}),$$
and if the corresponding $\lambda$-structure over a ring is special,
then $$nB_{n}=\sum_{d|n}\mu(d)\Psi_d(C_{n\over d}),$$
wher $\mu$ is a Moebius function.
\end{theorem}

\begin{proof}
Remark that $$(1-t^{k})^{B_k}{d\over dt}(1-t^{k})^{-B_k}=kt^{k-1}\sum_{n=1}^{\infty}\Psi_{n}(B)t^{k(n-1)}=\sum_{n=1}^{\infty}\Psi_n(kB_k)t^{kn-1},$$
so the equation (\ref{prod})
is equivalent to the equation
$$\sum_{n=1}^{\infty}C_{n}t^{n-1}=\sum_{k,m=1}^{\infty}\Psi_m(kB_k)t^{km-1},$$
that is $$\sum_{km=n}\Psi_{m}(kB_k)=C_n.$$
We have $B_1=C_1, 2B_2+\Psi_2(B_1)=C_2$ etc, so the solution for this system of equations is unique.
From the other hand, it is easy to see that the expressions for $B_k$ in the statement of the lemma
satisfies the last equation.
\end{proof}

{\bf Example 1.} Let $A(t)=\exp(at)$, then $C(t)=a,$
so $C_1=a,$ $C_i=0,i>1.$
Suppose that the $\lambda$-structure is special.
Therefore $$nB_n=\sum_{d|n}\mu(d)\Psi_d(C_{n\over d})=\mu(n)\Psi_{n}(a), B_n={\mu(n)\over n}\Psi_{n}(a).$$ Hence,
$$\exp(at)=\prod_{n=1}^{\infty}(1-t^n)^{-{\mu(n)\over n}\Psi_{n}(a)},$$
where right hand side is considered in a sense of the power structure.
For example, for $a=1$ we get the equality of formal power series
$$e^t=\prod_{n=1}^{\infty}(1-t^n)^{-{\mu(n)\over n}}.$$

{\bf Example 1a.} Analogously to the previous example one can also prove a couple of curious identities

$$\prod_{k=1}^{\infty}(1-t^{k})^{-{\varphi(k)\over k}}=e^{t\over 1-t},$$

where $\varphi(n)$ is the Euler function, that is the number of integers less than $n$ and coprime with $n$,
and

$$\prod_{g.c.d(k,m)=1}(1-x^{k}y^{m})^{1\over k}=(1-x)^{y\over 1-y}.$$

{\bf Example 2.} Let $A(t)=1+at,$ then $C(t)={a\over 1+at},$
so $C_j=-(-a)^{j}.$ Suppose that the $\lambda$-structure is special, then
$$nB_{n}=-\sum_{d|n}\mu(d)\Psi_d((-a)^{n\over d}).$$
Therefore $${d\over dt}\ln((1+at)^x)=\sum_{n=1}^{\infty}D_{n}t^{n-1},$$
where $$D_n=\sum_{km=n}\Psi_{m}(kB_k\cdot x)=-\sum_{km=n}\Psi_{m}(\sum_{d|k}\mu(d)\Psi_d((-a)^{k\over d})\cdot x)=$$
$$-\sum_{dsm=n}\mu(d)\Psi_{md}((-a)^s)\Psi_{m}(x).$$
Hence,
$$\ln((1+at)^x)=-\sum_{d,s,m=1}^{\infty}\mu(d)\Psi_{md}((-a)^s)\Psi_{m}(x){t^{dsm}\over dsm}=$$
$$\sum_{d,m=1}^{\infty}{\mu(d)\over dm}\ln(1+\Psi_{md}(a)t^{md})\Psi_m(x)=
\sum_{n=1}^{\infty}\ln(1+\Psi_{n}(a)t^{n}){1\over n}\sum_{m|n}
\mu({n\over m})\Psi_m(x).$$
We get a formula
$$(1+at)^x=\prod_{n=1}^{\infty}(1+\Psi_{n}(a)t^{n})^{{1\over n}\sum_{m|n}
\mu({n\over m})\Psi_m(x)},$$
where powers in the left hand side are considered in a sense of the power structure, and in the right hand side -- in a "ususal sense of the exponent of the logarithm".

\bigskip

Let's consider, for example, a polynomial ring $\Lambda$  of the infinite number of variables $x_1, x_2,\ldots$ with the integer coefficients. Let $p_k=x_1^k+x_2^k+\ldots$ denote the Newton symmetric polynomials.
Then $\Psi_k(p_1)=p_k.$ 

Let $X$ be a quasiprojective variety,
let $F(X,n)$ denote a set of ordered $n$-tuples of points on $X$, $e^{S_n}_{F(X,n)}(u,v)$ is an equivariant Hodge-Deligne polynomial (\cite{getzler}) for the natural $S_n$-action on $F(X, n)$. Let $B(X,n)$ be a set of unordered tuples of distinct points on $X$.
E. Getzler (\cite{getzler}) proved that the following identity is true
$$1+\sum_{n=1}^{\infty}t^{n}e^{S_n}_{F(X,n)}(u,v)=\prod_{n=1}^{\infty}(1+p_nt^{n})^{{1\over n}\sum_{m|n}
\mu({n\over m})\Psi_m(e_X(u,v))}.$$

Using the discussion above, we can rewrite this equation in a simpler form
\begin{equation}
\label{conf}
1+\sum_{n=1}^{\infty}t^{n}e^{S_n}_{F(X,n)}(u,v)=(1+p_1t)^{e_X(u,v)},
\end{equation}
where the right hand side is considered in a sense of the power structure over a ring
$\Lambda\otimes\mathbb{Z}(u,v)$.

{\bf Example 3a.} One can prove (\cite{macd}) that after the change of all $p_i$ to 1 in the character of a representation he gets a (virtual) multiplicity of the trivial representation in a given one.  Thus if we change all  $p_i$ to 1 in the Getzler's formula,
we'll get a generating function for the dimensions of the $S_n$-invariant subspaces in the cohomologies of $F(X,n)$,
that is a generating function for the Poincare polynomials (with the compact support) of the quotients
$B(X,n)=F(X,n)/S_n$. Therefore
$$1+\sum_{n=1}^{\infty}t^nP_{B(X,n)}(q)=
\prod_{n=1}^{\infty}(1+t^{n})^{{1\over n}\sum_{m|n}\mu({n\over m})\Psi_m(P(X))},$$
what coincides with $(1+t)^{P(X)}$ in a sense of the power structure. This coincidence is not by chance, since from the geometric interpretation of the power structure over the Grothendieck ring of varieties
$$1+\sum_{n=1}^{\infty}t^n[B(X,n)]=(1+t)^{[X]}.$$

If $b_k$ is a $k$th Betti number of the variety $X$,
then
$$(1+t)^{P(X)}={(1-t^2)^{P(X)}\over (1-t)^{P(X)}}=\prod_{k=0}^{\infty}{(1-t^2q^k)^{(-1)^k b_k}\over (1-tq^k)^{(-1)^k b_k}}.$$

{\bf Example 3b.} A loop on the $B(X,n)$ corresponds to an automorphism of the covering  $F(X,n)\rightarrow B(X,n)$, that is an element of the symmetric group
$S_n$. Therefore, for example, the sign representation of $S_n$ corresponds to some one-dimensional representation of 
$\pi_1(B(X,n))$.

It is known (\cite{macd}) that the change of $p_i$ to $(-1)^{i-1}$ in the character of a representation gives a (virtual) multiplicity of the sign representation in a given one. Hence if we change all $p_i$ to $(-1)^{i-1}$ in the Getzler's formula, we'll get a generating function for the Poincare polynomials (with the compact support) of
$B(X,n)$ with the coefficients in the sign representation. Therefore
$$1+\sum_{n=1}^{\infty}t^n\sum_{k=0}^{\infty}(-1)^{k}q^k\dim H^{k}_{c}(B(X,n),\pm \mathbb{C})=
\prod_{n=1}^{\infty}(1+(-1)^{n-1}t^{n})^{{1\over n}\sum_{m|n}\mu(n/m)\Psi_m(P(X))},$$
what coincides with $(1-u)^{P(X)}|_{u=-t}$ in a sense of the power structure. If $b_k$ is a $k$th Betti number of $X$,
then
$$(1-u)^{P(X)}|_{u=-t}=\prod_{k=0}^{\infty}(1-uq^k)^{(-1)^k b_k}|_{u=-t}=\prod_{k=0}^{\infty}(1+tq^k)^{(-1)^k b_k}.$$

{\bf Example 3c.}  One can prove (\cite{macd}) that the coefficient at $p_1^n$ in the character of a $S_n$-representation equals to the (virtual) dimension of this representation, multiplied by $n!$. Therefore if we change all  $p_i$ to 0 for $i>1$ in the Getzler's formula, we obtain the exponential generating function for the Poincare polynomials (with the compact support) of
$F(X,n)$. Therefore
$$1+\sum_{n=1}^{\infty}{t^n\over n!}P_{F(X,n)}(q)=(1+t)^{P(X)},$$
where the power is considered in the "usual" sense.

{\bf Example 4.} O. Tommasi (\cite{tommasi},\cite{bertom}) proved that the homologies of the moduli space of hyperelliptic curves of arbitrary genus are trivial (but with nonzero weight). Using the power srtucture, it is easy to check that the Hodge-Deligne polynomial of the moduli space of hypoerelliptic curves of genus 
 $g$ is equal to $(uv)^{4g-2}.$ 
 
 Every hyperelliptic curve is in 1-to-1 correspondence with an
 unordered $(2g+2)$-tuple of distinct points on $\mathbb{P}^1$ up to the action of the group $PGL(2,\mathbb{C})$.
Projectivisation  of the space of $2\times 2$-matrices is  isomorphic to $\mathbb{CP}^3$, its class in the Grothendieck ring equals to $\mathbb{L}^3+\mathbb{L}^2+\mathbb{L}+1.$
Degenerate matrices lay on the Segre quadric, whose class is equal to $(1+\mathbb{L})^2=1+2\mathbb{L}+\mathbb{L}^2.$
Hence,
$$[PGL(2,\mathbb{C})]=(\mathbb{L}^3+\mathbb{L}^2+\mathbb{L}+1)-(1+2\mathbb{L}+\mathbb{L}^2)=(\mathbb{L}^3-\mathbb{L}).$$
From the other hand, the generating function for the classes of the unordered tuples of points on $\mathbb{P}^1$ has a following form:
  $$(1+t)^{1+\mathbb{L}}=(1+t)\cdot{(1-t)^{-\mathbb{L}}\over (1-t^2)^{-\mathbb{L}}}={(1-\mathbb{L}t^2)(1+t)\over (1-\mathbb{L}t)}.$$
  The coefficient at $t^k$  in the decomoposition of this function equals to $\mathbb{L}^k-\mathbb{L}^{k-2}$ for $k\ge 4$.

Therefore the class of the moduli space of hyperelliptic curves in the Grothendieck ring equals to  
$${\mathbb{L}^{2g+2}-\mathbb{L}^{2g}\over \mathbb{L}^3-\mathbb{L}}=\mathbb{L}^{2g-1},$$
so we get the desired statement about the Hodge-deligne polynomial.

\bigskip

\bigskip

It turns out that the equation from the example 2 can be generalized to the arbitrary series.

\begin{theorem}
\begin{equation}
\label{powprod}
(1+A_1t+A_2t^2+\ldots)^X=\prod_{n=1}^{\infty}(1+\Psi_n(A_1)t^n+\Psi_n(A_2)t^{2n}+\ldots)^{{1\over n}\sum_{m|n}
\mu({n\over m})\Psi_m(X)},
\end{equation}
where powers in the left hand side is considered in a sense of a power structure, and in the right hand side -- in a sense of the "exponent of the logarithm".
\end{theorem}

\begin{proof}
Let us prove first that the right hand side of the equation (\ref{powprod})
defines a power structure. The properties 1,3, 4,6,7 are obvious.
The property 2 follows from the equation $\sum_{m|n}
\mu(n/m)=\delta_{n,1}.$ Let us prove the property 5.
From the viepoint of this (conjectural) structure
$$(((A(t))^X)^Y=[\prod_{n=1}^{\infty}(1+\Psi_n(A_1)t^n+\Psi_n(A_2)t^{2n}+\ldots)^{{1\over n}\sum_{m|n}
\mu(n/m)\Psi_m(X)}]^Y=$$ $$=\prod_{k=1}^{\infty}\prod_{n=1}^{\infty}
\Psi_{k}[(\Psi_{n}(A)(t^{kn}))^{{1\over n}\sum_{m|n}\mu(n/m)\Psi_m(X)}]^{{1\over k}\sum_{l|k}\mu(k/l)\Psi_l(Y)}=$$
$$=\prod_{k=1}^{\infty}\prod_{n=1}^{\infty}(\Psi_{kn}(A)(t^{kn}))^{{1\over kn}\sum_{m|n}\sum_{l|k}\mu(n/m)\mu(k/l)\Psi_km(X)\Psi_l(Y)}.$$
Let $a=kn, b=km$.Then $n/m=a/b$.
Let us note that $$\sum_{k: k|b, k\vdots l}\mu(k/l)=\sum_{f|(b/l)}\mu(f)=\delta_{b/l,1}.$$
Therefore $$\sum_{kn=a}\sum_{m|n}\sum_{l|k}\mu(n/m)\mu(k/l)\Psi_{km}(X)\Psi_l(Y)=\sum_{b|a}\mu(a/b)\Psi_{b}(X)\Psi_b(Y)=\sum_{b|a}\mu(a/b)\Psi_b(XY).$$
Hence $$\prod_{k=1}^{\infty}\prod_{n=1}^{\infty}(\Psi_{kn}(A)(t^{kn}))^{{1\over kn}\sum_{m|n}\sum_{l|k}\mu(n/m)\mu(k/l)\Psi_km(X)\Psi_l(Y)}=
\prod_{a=1}^{\infty}(\Psi_{a}(A)(t^{a})^{{1\over a}\sum_{b|a}\mu(a/b)\Psi_b(XY)}.$$

Therefore the property  5 is also satisfied, and the right hand side of the equation (\ref{powprod}) defines a power structure. What rests to prove is the coincidence of this structure with the initial one at $A(t)=(1-t)^{-X}$.
Logarithming the expression
\begin{equation}
\label{fla}
\prod_{n=1}^{\infty}(1-t^n)^{-{1\over n}\sum_{m|n}
\mu(n/m)\Psi_m(X)},
\end{equation}
we get
$$-\sum_{n=1}^{\infty}{1\over n}\sum_{m|n}
\mu(n/m)\Psi_m(X)\ln(1-t^n)=\sum_{n=1}^{\infty}\sum_{m|n}\sum_{k=1}^{\infty}{1\over n}\mu(n/m){1\over k}t^{kn}\Psi_{m}(X)=$$
$$=\sum_{s,m=1}^{\infty}\sum_{l|s/m}\mu(l){1\over s}t^{s}\Psi_{m}(X)=
\sum_{s=1}^{\infty}{1\over s}t^{s}\Psi_{m}(X).$$
Therefore the logarithmic derivative of (\ref{fla}) equals to
$$\sum_{s=1}^{\infty}t^{s-1}\Psi_{m}(X),$$
what finishes the proof.
\end{proof}

{\bf Example 5.}
Let $P_N$ denote the projectivization of the set of polynomials of degree
$N$ of $n$ variables, and let $Irr_N$ be the projectivization of the set of the irreducible polynomials of degree $N$ of $n$ variables.

The class of $P_N$ in the Grothendieck ring of varieties equals to
$$[P_N]={\mathbb{L}^{n\choose n+N}-\mathbb{L}^{n\choose n+N-1}\over \mathbb{L}-1}.$$
Let $$P(\mathbb{L},t)=1+\sum_{N=1}^{\infty}[P_N]t^N.$$

Let $$P(\mathbb{L},t)^{-1}{\partial\over \partial t}P(\mathbb{L},t)=\sum_{n=1}^{\infty}C_{n}(\mathbb{L})t^{n-1}.$$

\begin{theorem}
$$n[Irr_n]=\sum_{d|n}\mu(d)C_{n\over d}(\mathbb{L}^{d}).$$
\end{theorem}

\begin{proof}
Since in the polynomial ring the decomposition into the irreducible factors is unique up to multiplication by a constant,
one can define a set $P_{i_1i_2\ldots}$ of polynomials which are products
of $i_1$ irreducible factors of degree 1
1, $i_2$ of degree 2 ($i_k=0$ for $k$ large enough).

Note that $$P_{i_1i_2\ldots}=(S^{i_1}Irr_1)\times(S^{i_2}Irr_2)\times\ldots$$
and $$P_N=\sqcup_{i_1+2i_2+\ldots=N}P_{i_1i_2\ldots},$$
so $$P(\mathbb{L},t)=\sum_{N=0}^{\infty}[P_N]t^N=\prod_{k=1}^{\infty}(1-t^k)^{-[Irr_k]}.$$

Now the proposition of the theorem follows from the lemma 1 and the example before it.
\end{proof}

\begin{corollary}
The Hodge-Deligne polynomial of the projectivization of the set of irreducible polynomials is defined by a formula
$$n\cdot e_{Irr_n}(u,v)=\sum_{d|n}\mu(d)C_{n\over d}(u^{d}v^{d}).$$
\end{corollary}

One can calculate the Euler characteristic of  $Irr_i$. 
$$\prod_{i=1}^{\infty}(1-t^i)^{-\chi(Irr_i)}={1\over (1-t)^n}.$$
Taking the logarithms, $$\sum_{i}\chi(Irr_i)\ln(1-t^i)=n\ln(1-t),$$
and $$\sum_{n}\chi(Irr_i)\sum_{l}t^{il}/l=n\sum_{k}t^k/k.$$
Comparing the coefficients, we get
$$\sum_{i|k}i\chi(Irr_i)=n,$$
So by the Moebius inversion formula
$$\chi(Irr_i)={2\over i}\sum_{k|i}\mu(i),$$
what equals to $n$ for $i=1$ and to 0 for $i>1$.

\section{Plethysms and representations of the symmetric groups.}

Let $\Lambda$ be the ring of symmetric polynomials of the infinite number of variables. Over $\Lambda$, as over any polynomial ring, there is a natural  $\lambda$-structure.
It is easy to check that this structure is special.

For any $\lambda$-ring $R$ one can construct a natural map, sending a pair of elements $f\in\Lambda$ and $X\in R$ to an element
$f\circ X\in R$ such that the following properties are satisfied:

1)$(f_1+f_2)\circ X=f_1\circ X+f_2\circ X;$

2)$(f_1f_2)\circ X=(f_1\circ X)(f_2\circ X);$

3)$p_k\circ X=\Psi_{k}(X).$

An easy check shows that
 $$(1-t)^{-X}=\sum_{k=0}^{\infty}h_{k}\circ X t^k.$$

One can prove that a $\lambda$-ring is special if and only if for any $f,g, X$ $(f\circ g)\circ X=f\circ(g\circ X).$

Consider a direct sum of representation rings of groups  $S_k$ over all $k$. For each representation $V$ of the group $S_k$ one can construct its character -- it is a homogeneous polynomial from  $\Lambda$ of degree $k$, which can be defined by a formula
$$ch(V)=\sum_{\sigma\in S_{k}}Tr(\sigma|_{V})\cdot p_1^{i_1(\sigma)}\cdot\ldots\cdot p_k^{i_k(\sigma)},$$
where $i_s(\sigma)$ -- is a number of cycles of length $s$ in a representation
$\sigma$.

It is proved in \cite{macd}, that the  "natural" operation for the representations corresponds to the natural operations for their characters: the character of the direct sum of representations is a sum of their characters, 
if $V$ is a representation of $S_k$, and $W$ is a representation of $S_m$,
then the character of the representation $Ind_{S_{k}\times S_{m}}^{S_{k+m}}(V\otimes W)$ is equal to a product of characters of $V$ and $W$. Moreover, on the representations there is a so-called plethysm operation: if $V$ is a representation of $S_k$, and $W$ is a representation of  $S_m$, then on the product $V\otimes W^{\otimes k}$ 
there is a natural action of the semidirect product of the groups $S_k$ and $(S_{m})^{k}$.
The plethysm is a representation $V\circ W=Ind_{S_k\propto (S_{m})^{k}}^{S_{km}}(V\otimes W^{\otimes k})$. It turns out (\cite{macd}),
that $$ch(V\circ W)=ch(V)\circ ch(W).$$ In particular, if
$V_k$ is a trivial one-dimensional representation of the $S_k$, then $ch(V_k)=h_{k}$,
and $ch(V_k\circ W)=h_k\circ ch(W).$ The generating function for characters of such representations is equal to 
$$\sum_{k=0}^{\infty}h_k\circ ch(W)=(1-t)^{-ch(W)}.$$

Consider a direct sum of the Grothendieck rings of the representations of all symmeric groups: $$R=\oplus_{n=1}^{\infty}R(S_n).$$
The multiplication on the ring $R$ is given by the formula
$$V\odot W=Ind_{S_{k}\times S_{m}}^{S_{k+m}}(V\otimes W),$$
if $V$ is a representation of  $S_k$, and $W$ is a representation of $S_m$
(the product is a representation of $S_{k+m}$).

It is proved in \cite{macd},\cite{knut}, that the character map is an isomorphism between this ring and the ring $\Lambda$.
From the above discussion it follows that the power structure over the polynomial ring $\Lambda$ corresponds to a structure over the ring $R$, such that
$$(1-t)^{-W}=\sum_{k=0}^{\infty}t^{k}\cdot Ind_{S_k\propto (S_{m})^{k}}^{S_{km}}W^{\otimes k}$$
for a representation  $W$ of the group $S_m$.

\section{Moduli spaces of curves}

Since every complex curve of genus 2 is hyperelliptic, one can try to compute the $S_n$-equivariant Euler characteristic of the moduli space of genus 2 curves with $n$ marked points. 

Consider the forgetful map 
$\mathcal{M}_{2,n}\rightarrow \mathcal{M}_2$. It is not a locally trivial fibration, since a curve can have a nontrivial automorphism group (for example, every hyperelliptic curve has a notrivial automorphism -- a hyperelliptic involution), and a true fiber of the forgetful map is a quotient of the space of distinct unordered points on a curve by the action of an automorphism group of a curve.

Consider a following problem: let a finite group $G$ acts on a variety  $X$.  Let $X_k(g)$ denote the set of points with orbit of length $k$ for the action of an element $g\in G$.

\begin{theorem}
$$\sum_{n=0}^{\infty}t^{n}\chi^{S_n}(F(X,n)/G)={1\over |G|}\sum_{g\in G}\prod_{k}(1+p_kt^{k})^{\chi(X_k(g))\over k}.$$
\end{theorem}

\begin{proof}
Let $R\in Rep(G)$ be an alternating sum of the cohomologies of $X$ as representations of $G$ ($R$ belongs to the Grothendieck ring of the representations of the group $G$). By the Getzler's formula the analogous sum for the $S_n$-equivariant cohomologies of $F(X,n)$ equals to
$$\prod_{k=1}^{\infty}(1+p_kt^{k})^{{1\over k}\sum_{d|k}\mu(k/d)\Psi_d(R)},$$
and its character at an element $g$ equals to
$$\xi(g)=\prod_{k=1}^{\infty}(1+p_kt^{k})^{{1\over k}\sum_{d|k}\mu(k/d)\chi(g^d)},$$
if $\chi$ is a character of a representation $R$. Note that from the Lefschetz theorem it follows that
$\chi(e)=\chi(X)$, and for all other $g$ $\chi(g)$ is equal to the Euler characteristic of the fixed point set of $g$.
The dimension of a $G$-invariant part in the cohomologies of  $F(x,n)$ equals to
$${1\over |G|}\sum_{g\in G}\xi(g)$$.

The number $\chi(g^d)$ equals to the Euler characteristic of the $g^d$-fixed point set, that is 
$\sum_{l|d}\chi(X_l(g)).$
Therefore, $${1\over k}\sum_{d|k}\mu(k/d)\chi(g^d)=
{1\over k}\sum_{d|k}\mu(k/d)\sum_{l|d}\chi(X_l(g))={1\over k}\sum_{l|k}\chi(X_l(g))\sum_{d|k}\mu(d)=
{1\over k}\chi(X_k(g)).$$
Объединяя эти ответы, получим
$$\xi(g)=\prod_{k}(1+p_kt^{k})^{\chi(X_k(g))\over k},$$
what finishes the proof.
\end{proof}

\begin{corollary}
$$\sum_{n=0}^{\infty}{t^{n}\over n!}\chi(F(X,n)/G)={1\over |G|}((1+t)^{\chi(X)}+\sum_{g\neq e}(1+t)^{L(g)}),$$
where $L(g)=\chi(X_1(g))$ is an Euler characteristic of the fixed point set of $g$.
\end{corollary}

To calculate the equivariant Euler characteristic of $\mathcal{M}_{2,n}$ one has
to describe hyperelliptic curves with additional symmetries. For this one has to describe all  6-tuples of points on
$\mathbb{CP}^1$, having nontrivial symmetry groups and to compute the Euler characteristics of the moduli spaces of the corresponding curves.

The sum of these Euler characteristics equals to 1 since the theorem of O. Tommasi
(\cite{tommasi},\cite{bertom}) states that the homologies of the moduli space of hyperelliptic curves are trivial.

The sum of the Euler characteristics of strata, divided by the orders of corresponding symmetry groups (the orbifold Euler characteristic), is equal to
${-1\over 240}$
in agreemeent with the Harer-Zagier formula (\cite{harzag}):
$$\chi_{orb}(\mathcal{M}_{g,n})=(-1)^{n}{(2g-3+n)!(2g-1)\over (2g)!}B_{2g},$$
where $B_{2g}$ are Bernoulli numbers. For $g=2$ $B_{2g}={-1\over 30},$
so  $$\chi_{orb}(\mathcal{M}_{2,0})={1\cdot 3\over 4!}\cdot {-1\over 30}={-1\over 240}.$$

Summing up the answers for different strata, we get (\cite{moduli}) the following

\begin{proposition}
$$\sum_{n=0}^{\infty}t^n\chi^{S_n}(\mathcal{M}_{2,n}))={-1\over 240}(1+p_1t)^{-2}-{1\over 240}(1+p_1t)^6(1+p_2t^2)^{-4}+$$
$$+{2\over 5}(1+p_1t)^3(1+p_5t^5)^{-1}+{2\over 5}(1+p_1t)(1+p_2t^2)(1+p_5t^5)(1+p_{10}t^{10})^{-1}+$$
$$+{1\over 6}(1+p_1t)^2(1+p_2t^2)(1+p_6t^6)^{-1}-{1\over 12}(1+p_1t)^4(1+p_3t^3)^{-2}-$$
$$-{1\over 12}(1+p_2t^2)^2(1+p_3t^3)^2(1+p_6t^6)^{-2}+{1\over 12}(1+p_1t)^2(1+p_2t^2)^{-2}+$$
$$+{1\over 4}(1+p_1t)^2(1+p_4t^4)(1+p_8t^8)^{-1}-{1\over 8}(1+p_1t)^2(1+p_2t^2)^2(1+p_4t^4)^{-2}.$$
\end{proposition}

Up to 4 points we get

$$\sum_{n=0}^{\infty}t^n\chi^{S_n}(\mathcal{M}_{2,n}))=
1+2p_1\cdot t+p_1^2\cdot t^2+({1\over 2}p_4+{2\over 3}p_1p_3-{1\over 6}p_1^4)\cdot t^4+\ldots=$$
$$1+2s_1\cdot t^4+(s_1+s_2)\cdot t^2+(s_4-s_{3,1}-s_{2,2})\cdot t^4+\ldots,$$

what coincides with the answers obtained by Getzler (\cite{TRR}), and Tommasi and Bergstrom (\cite{bertom}).

\vspace{1cm}

Author is grateful to S. Gusein-Zade, M. Kazaryan and S. Lando for  lots of useful discussions.

\renewcommand{\refname}{Bibliography}

Moscow State University,\newline
Department of Mathematics and Mechanics.

Е.mail: gorsky@mccme.ru.

\end{document}